\documentclass[12]{article}
\usepackage{amsmath}
\usepackage{amsfonts}
\usepackage{amssymb}
\usepackage{amsthm}
\usepackage{color}

\def\binom#1#2{{#1}\choose{#2}}
\newcommand{\R}{{\rm I}\kern-0.18em{\rm R}}
\newcommand{\1}{{\rm 1}\kern-0.25em{\rm I}}
\newcommand{\E}{{\rm I}\kern-0.18em{\rm E}}
\newcommand{\p}{{\rm I}\kern-0.18em{\rm P}}
\makeatletter
\renewcommand{\@cite}[ 2]{{#1\if@tempswa , #2\fi}}
\renewcommand{\@biblabel}[1]{ \hfill}
\makeatother

\usepackage{epsfig}

\usepackage{fancyhdr}

\usepackage{graphics}

\usepackage{amsopn}  

\usepackage{amssymb}
\usepackage{amsxtra}
\usepackage{amsfonts}

\usepackage{afterpage}
\usepackage{eucal}
\usepackage{eufrak}

\usepackage{enumerate}

\allowbreak

\def\fnote#1{\footnote}

\newcommand{\bea}{\begin{eqnarray}}

\newcommand{\eea}{\end{eqnarray}}

\newcommand{\beas}{\begin{eqnarray*}}
\newcommand{\eeas}{\end{eqnarray*}}

\author{Lev B Klebanov and Lenka Slamova}
\title{Discrete Stable and Casual Stable Random Variables}
\date{Stability of Stochastic Models, 2014}
\begin{document}

\maketitle

\begin{abstract}
Here we introduce some new classes of discrete stable random variables, which are useful for understanding of a new general notion of stability of random variables called us as casual stability. There are given some examples of casual and discrete stable random variables. We also propose a class of discrete stable random variables for a description of rating of scientific work.
\end{abstract}


\section{Classical Strictly Stable Random Variables}
Let $X_1, \ldots , X_n, \ldots$ be  i.i.d. r.v.  $ \forall n \in N$ there is $a_n \in (0,1)$ such that
\[ X_1 \stackrel{d}{=} a_n \sum_{k=1}^n X_k. \]
In this case we say that $X_1$ is strictly stable random variable (r.v).

There are many applications to physics, astronomy, finance. However, the definition is not applicable to r.v.s taking positive integer values. The notion of discrete stability for lattice random variables on nonnegative integers was introduced in Steutel and van Harn (1979).

\section{Discrete Strictly Stable Random Variables}
Let $X_1, \ldots , X_n, \ldots$ be  independent and identically distributed (i.i.d.) r.v. taking the values in $N$. Let $\varepsilon_{i,j}(p)$  be i.i.d. r.v.s with Bernoully distribution, taking values $0$ with probability $1-p$ and $1$ with probability $p$. Define
\[ \tilde{ X_j}(p) =\sum_{i=1}^{X_j}\varepsilon_{i,j}(p) \]
as a normalization of r.v. $X_j$. This normalization can be considered in the following way:
$X_j = 1+\ldots +1$ $ X_j$ times, so we have $X_j$ units (particles) passing through some stuff. Each particle  may be absorbed by this stuff with probability $1-p$, and not absorbed with probability $p$. The number of non-absorbed particles is r.v. $\tilde{X_j}(p)$.

 Steutel and van Harn (1979) say that the r.v. $X_1$ is a (strictly) discrete stable r.v. if $\forall$ $n \in N$ there is $p(n) \in (0,1)$ such that
 \[ X_1 \stackrel{d}{=}\sum_{j=1}^{n}\tilde{X_j}(p_n). \]
 The p.g.f. of discrete stable r.v. satisfies to the equation
 \[ \mathcal{P}(z)=\mathcal{P}((1-p(n)+p(n)z)^{n} \]
  and has the following form
 \[ \mathcal{P}(z) = \exp \{- \lambda (1-z)^\alpha\}, \]
 where $\lambda >0$, and $\alpha \in (0,1]$. Because this distribution is concentrated on positive semi-axis, the parameter $\alpha$ cannot be greater than $1$.

\section{Another definition of discrete stability}
However, we think, that such normalization procedure is not unique possible. Of course, it seems natural, that some of particles passing through special stuff may generate other particles. The fact that the stuff is absorbing may be expressed by the fact that the mean number of particles on the other side of the stuff is less that $1$. Mathematically, we may change Bernoulli distribution of $\varepsilon$\rq{}s by another distribution, with mean value less that $1$. Suppose that $\varepsilon$\rq{}s have p.g.f. $Q_n(z)$.  The p.g.f. of $X_1$ has to satisfy the following equation:
\[ \mathcal{P}(z) =\mathcal{P}(Q_n(z))^n, \; \forall n \in N. \]
It is clear, that $P$ is probability generating function (p.g.f.) of an infinite divisible (i.d.) distribution. Semigroup generated by the family $\{Q_n\}$ with superposition operation has to be commutative.

Let us give some examples of new discrete stable distributions.

{\bf Example 1}. Define
\[ Q(z)= \Bigl( \frac{(1-p)+(p-\kappa)z^m)}{(1-p\kappa)-\kappa (1-p)z^m}\Bigr)^{1/m}, \]
$p=p(n)$. It can be shown, that $Q(z)$ is p.g.f., and the family of $Q(z)$ as a function of $p$ is a commutative semigroup. The parameters in $Q$ are a) for $m=1:$ $0<= \kappa <1$, $0<p<1$, and 
b) for $m \in N$, $m>1:$ $0<p<\kappa <1$.  Corresponding p.g.f. of discrete stable $X_1$ has the form
\[ \mathcal{P}(z) = \exp\{-\lambda \Bigl( \frac{1-z^m}{1-\kappa z^m}\Bigr)^{\gamma}\}. \]
The case $\kappa=0$, $m=1$ leads to Steutel and van Harn definition.

{\bf Example 2} Define
\[ Q(z) = \frac{2(b+T_p(\frac{(1+b)z-2b}{2-(1+b)z}))}{(1+b)(1+T_p(\frac{(1+b)z-2b}{2-(1+b)z}))},\]
where $p \in (0,1)$, $b \in (-1,1)$, and $T_p(x)=\cos (p \arccos x)$ is for $p=1/n$ a function inverse to Chebyshev polynomial; $p=p(n)$. The function $Q(z)$ is a p.g.f. (it was difficult to prove for us), and corresponding semigroup is commutative. P.g.f. of discrete stable $X_1$  has the form:
\[ \mathcal {P}(z)= \exp\{- \lambda (\arccos \frac{(1+b)z-2b}{2-(1+b)z})^{\gamma}\},\]
where $\lambda >0$, $\gamma \in (0,2]$ and $b \in (-1,1)$.

\section{Rating  of  scientific  work  and  discrete  stable  distributions}
It  is  very  often  to base  rating  of  scientific  work  on  the  number  of  citations  of  corresponding  paper,  author  or  journal, in which  the  paper  was published  (so - called  impact factor).  Below  we  give  a  model  of  the  distribution of citations  number, and  show  a  connection  with  discrete  stable  distributions.

At  start  let  us  consider  a simplest  model  of  paper  publication.    We  are  consider  the case  only when  there  is at  least  one  publication  at start (in opposite  case  there  will  be  no citations  at all).   Let  the  probability  of  a rejection of  a paper  is  q. Then,  the  probability  to have  exactly  $k$ published  papers is  the  probability  of $k - 1$ acceptions  (it  is  $(1 - q)^(k - 1)$), (we  suppose  one  paper  was  published)  and one rejection (with  probability  $q$), so, this  probability  is  $q (1 -q)^(k - 1)$. In other  words,  we  have  geometric  distribution  for  the  number  of  published  papers.  It  has p.g.f. 
\[ Q(z) =\frac{q z}{1 - (1 - q) z}. \]

Suppose  now, that each published paper generates some citations. It is more or less clear, that the probability for paper to be cited depends on the number of previous its citations. Suppose, that the probability that the paper having $k -1$ citations will be not cited again is $p/k$, where $p$ is the probability the paper will be not cited at all. Therefore the probability that the paper will cited exactly $k$ times  is 
\[ \frac{p}{k} \prod_{j=1}^{k-1}(1-p/j) = p (1-p)_{k-1}/k! ,\]
where $(a)_n=a(a+1)\cdots (a+n-1)$ is the Pochhammer symbol. Therefore,  the  p.g.f. of the distribution of citations of one paper is
\[  \sum_{k=1}^{\infty} \frac{p (1-p)_{k-1}}{k!}z^k = 1-(1-z)^p. \]
This is well - known Sibuya distribution with parameter $p$. So, the number of citations of one publication has Sibuya distribution with parameter $p$.

The  p.g.f. of the number of citations for all paper coming from one author is the superposition of p.g.f.  
of Sibuya distribution with p.g.f. of geometric distribution, that is 
\[ 1-\Bigl( 1- \frac{qz}{1-(1-q)z} \Bigr)^{p}\]
with parameters $q \in (0, 1]$  and $p \in (0, 1]$. So, we have p.g.f of the number of citations of one author. It is easy to see, that corresponding distribution has a heavy tail (its limit behavior is of order $1/k^p$ as $k \to \infty $. Suppose now that we are interested in the distribution of the number of citations  in some field of science. It is natural to suppose, that the number of scientists having publications in this field has Poisson distribution with parameter $\lambda$. Then the p.g.f. of the number of all publications in  the  field is superposition of Poissonian  p.g.f. with p.g.f. described above, that is
\[ \exp \left ( -\lambda \Bigl( \frac{1-z}{1-(1-q)z} \Bigr)^{p}\right). \]
It  is  p.g.f. of discrete stable distribution with parameters $\lambda >0$, $p \in (0, 1]$ and $\kappa =1-q$.

We  can  easily  see, that  the  distribution  does  not have  finite  first  moment  if  $0<p< 1$.  It  has mode at zero, and finite median. So, on any empirical data we will see, that the empirical mean is much larger than empirical median. Also, many citations will respond to a (relatively) small number of publications, while the main part of publications will have small number of citations. This big difference between scientists is explained in our model just by random nature of publication and citation processes (of course, personal differences are also included in this randomness). Therefore, the ranking of scientists, scientific institutions or journals may not be based on the citation number. Such ranking will often produce random mistakes.

\section{A general definition of stability for additive system}
The latest definitions of discrete stability lead us to a new definition of stability in (more or less) general case of ``additive system\rq\rq{}. Let us start with a definition for non-negative r.v. Suppose that $X$ is non-negative r.v. with c.d.f $F(x)$. Its Laplace transform has the form
\[ L(s)=\int_0^\infty \exp (-sx) dF(x), \]
or (what is the same)
\[ L(s)=\int_0^\infty (\exp (-s))^x dF(x). \]
The function $\exp (-s)$ is a Laplave transform of degenerate distribution concentrated at point $1$.

Similar to discrete case, we change the function  $\exp (-s)$ by another function which is a Laplace transform of a distribution, concentrated on positive semi-axis, say $g(s)$. We would like to have the definition of new type of stability (call it casual stability) as: for $\forall n \in N$ there is a Laplace transform $g_n(s)$ such that
\[L^n(-\log g_n(s)) =L(s). \]
This equation is restrictive. Form it we see, that corresponding r.v. has to be i.d., and Laplace transform $g_n$ also cannot be arbitrary. However, we can substitute under $L$ sign Laplace transform of any i.d. distribution. If cumulative distribution function (c.d.f.)  $F(x)$ is concentrated in positive integers, then we can substitute any Laplace transform $g_n$. Of course, not any possible substitution leads to casual stable distribution.

The sense of casual stability may be explained in the following way. Set us suppose that we have a (discrete or continuous) flow of some elements (``particles\rq\rq{}) passing through a stuff. Each element of the flow may generate some other elements or just disappear, according to a probability distribution with Laplace transform $g_n(s)$. The casual stability means that the number of elements (or their characteristics)  of $n$ such flows after passing through the stuff is equivalent (that is has the same distribution) as initial flow before passing this stuff. 

A little bit differently, we may say that an additive system in ``randomly similar\rq\rq{} to its initial element, so the system is ``randomly self similar\rq\rq{}.

 Let us describe some examples of positive casual-stable r.v.\rq{}s.
 
{\bf \it The first example} is given by choice of $g_n(s)=\exp(-a_n s)$. In this case we have just an ordinary normalization, and corresponding casual stable distributions coincide with ordinary positive stable distributions skewed to the right with the parameter $\alpha \in (0,1)$.

{\bf \it The second example} is connected to discrete stability. Namely, let us choose
\[ g(s)= \exp\left\{-\Bigl( \frac{(1-p)+(p-\kappa)\exp(-ms))}{(1-p\kappa)-\kappa (1-p)\exp(-ms}\Bigr)^{1/m}\right\}, \]
where $p=p(n)$, and we come to the {\bf Example 1} above with the same parameters. As before, we can obtain discrete stability in the sense by  Steutel and van Harn putting here $\kappa =0$ and $m=1$.

It is clear, how to construct a function $g(s)$ corresponding to {\bf Example 2}. We will not discuss this here.

{\bf \it The third example} is more interesting from our point of view. 

Let us consider positive stable r.v. with the index of stability $\alpha=1/m$, $m \in N$, move its Laplace transform on $h$ units to the right and make normalization on corresponding measure. We will have Laplace transform of corresponding tempered stable r.v. as:
\[ L(s) =\exp\{ -\lambda^{\alpha}(1+\tan \frac{\pi \alpha}{2})((s+h)^{\alpha}-h^{\alpha})\}. \]
The distribution with this Laplace transform appears to be causal stable with the function
\[ g_n(s) = \exp(h-(1/n (s+h)^{\alpha}+ (n-1)/n h^{\alpha})^{1/\alpha})\]
(it is possible to verify that $g_n$ is Laplace transform of a d.f. in the case when $1/\alpha \in N$). For $h \to 0$ we obtain a classical case of normalization: degenerated distribution at point $1/n^{1/\alpha}$. As a particular case we find that Inverse Gaussian distribution is casual stable too ($\alpha = 1/2$).

\section{r.v.\rq{}s of arbitrary sign}
It is possible to consider the case of some r.v.\rq{}s taking values on the whole real line. Let $X$ be a r.v. with c.d.f. $F(x)$. Ch.f. $f(t)$ of $X$ may be written as
\[  f(t)= \int_{-\infty}^{\infty}e^{i t x} dF(x)= \int_{0}^{\infty}(e^{-i t})^{|x|} d(1-F(-x))+ \int_{0}^{\infty}(e^{i t})^{x} dF(x).\]
Now we make ``random normalization\rq\rq{} in the following way. We change $e^{i t}$ in the second integral by a characteristic function $g(t)$, but it the first integral we change $e^{-i t}$ by $g(-t)$. The definition of casual stability in this case is obvious now. The verification of the fact that a ch.f. is casual stable is generally more difficult than for the case of positive r.v.\rq{}s. However, we can prove, for example, that the Laplace distribution is casual stable.  The same fact holds for Linnik distribution as well.

\section{Casual $\nu$-stable distributions}
If in classical definition of stable distributions use a random number $\nu_p$ of summands instead of $n$, we come to so-called $\nu$-stability, and, particularly (for geometrically distributed $\nu$), to geo-stable distributions. Of course, casual stability can be generalized for the case of a random number of random summands (see Klebanov L.B., Maniya G.M., Melamed I.A., and Klebanov L.B., Rachev S.T )  Namely, it is possibly to use for this aim an isomorphism between i.d. and $\nu$-i.d. distributions: ch. f.  $f(t)$ is $\nu$-i.d. ch.f. if and only if
 \[ f(t) = \varphi(-\log(h(t))), \]
 where $h(t)$ is classical i.d. ch.f., and $\varphi$ is standard solution of Poincare equation.  Really, from definition we see that any casual stable distribution is i.d., and we can apply mentioned isomorphism to any such distribution.

\section{Convergence to casual stable distributions}
Here we give a limit theorem for convergence to casual stable distribution in the case of positive r.v.\rq{}s. General case can be considered in similar way, but the formulations appears to be more complicated.

Suppose that $L(s)$ is Laplace transform of a positive r.v., which is casual stable with Laplace transformations $g_n(s)$, that is $L(s)=L^n(-Log\; g(s))$, $\forall n \in N$.

{\bf Theorem} {\it Suppose that $h(s)$ is Laplace transform such that $\sup_{s>0}|h(s)-L(s)|/s^a<\infty$ for some positive $a$. Suppose also that
\[  \sup_{s>0}\frac{n s^a}{|g_n^{-1}(\exp(-s))|^a} \to 0, \; n\to \infty.\]
Let $X_1, \ldots , X_n , \ldots $ be a sequence of i.i.d. r.v.\rq{}s with Laplace transform $h(s)$. Then
\[ \sum_{j=1}^{n}\tilde{X}_j(n) \stackrel{d}{\to} Y, \]
where $\tilde{X}_j(n) $ is a $g_n$-normalized r.v., and $Y$ is a r.v. with Laplace transform $L(s)$.}

Let us give an example of application of this Theorem. 

Consider a r.v. with gamma-distribution, which has Laplace transform
\[ L(s) = \frac{1}{(1+b s)^{\gamma}},\]
where parameters $b$ and $\gamma$ are positive. Show, that gamma-distributed r.v. $Y$ is casual stable. 
Really, if so, then from definition of casual stability we must have
\[ L^{n}(-\log \; g_n(s)) =L(s) \]
for any integer $n>1$. From here we find
\[ g_n(s) = \exp\left\{ \frac{1}{b}\Bigl((1-(1+b s)^{1/n}\Bigr)\right\}. \]
It is clear, that for any integer $n>1$ the function $g_n(s)$ is Laplace transform of some probability distribution. Therefore, $Y$ is casual stable r.v. 

Suppose now, that $h(s)$ is Laplace transform of a r.v. $X_1$ such that
\[ \sup_{s>0}\frac{|h(s)-L(s)|}{s^a} <\infty  \]
for some $a>1$. We have
\[  \sup_{s>0}\frac{n s^a}{|g_n^{-1}(\exp(-s))|^a} =\frac{n s^a b^a}{((1+b s)^n-1)^a} = \]
\[ =\left(  \frac{n^{1/a}}{\sum_{k=1}^{n}{\binom{n}{k}} b^{k-1}s^{k-1} } \right)^a  \leq \frac{1}{n^{a-1}} \to 0\]
as $n \to \infty$. All the conditions of previous Theorem are met. So, we may say that  
\[ \sum_{j=1}^{n}\tilde{X}_j(n) \stackrel{d}{\to} Y, \]
where r.v.\rq{}s $\tilde{X}_j(n)$ are i.i.d. with Laplace transform $h(-\log\;g_n(s))$.

\section*{Acknowledgments}
The authors were supported by the Grant P 203/12/0665  GACR.

\end{document}